%% file: v_final.tex
\documentclass[a4paper, 11pt]{amsart}
\usepackage{longtable, stackrel, multicol, fancybox, tcolorbox, stackrel}
\usepackage{bm,amsfonts,amsmath,amssymb,mathrsfs,bbm,amsthm}
\usepackage{graphicx, color,hyperref,eurosym, eucal,palatino}

\definecolor{darkblue}{rgb}{0.2,0.2,0.6}
\definecolor{darkblue2}{rgb}{0.2,0.2,0.9}
\definecolor{superdarkblue}{rgb}{0.2,0.2,0.3}
\hypersetup{colorlinks=true, linkcolor=darkblue,citecolor=darkblue,
    urlcolor = superdarkblue}
\definecolor{citegreen}{rgb}{0.2,0.7,0.2}
\usepackage{marginnote}
\usepackage{geometry}
\input commands.tex

\renewcommand{\aa}{\alpha}

\theoremstyle{definition}

\newcommand\soutD{\bgroup\markoverwith
{\textcolor{DarkGreen}{\rule[.5ex]{2pt}{1pt}}}\ULon}

\makeatletter
\newcommand{\vast}{\bBigg@{3}}
\newcommand{\Vast}{\bBigg@{5}}
\makeatother

\begin{document}
%
\title[Spectral optimization for Robin Laplacian]{Spectral optimization for Robin Laplacian \\ on domains admitting parallel coordinates}

\author{Pavel Exner}
\address{P. Exner: Doppler Institute for Mathematical Physics and Applied Mathematics, Czech Technical University, B\v rehov\'a 7, 11519 Prague, and Department of Theoretical Physics, Nuclear Physics Institute, Czech Academy of Sciences, 25068 \v{R}e\v{z}, Czechia}
\email{exner@ujf.cas.cz}

\author{Vladimir Lotoreichik}
\address{V. Lotoreichik: Department of Theoretical Physics, Nuclear Physics Institute, Czech Academy of Sciences, 25068 \v{R}e\v{z} Prague, Czechia}
\email{lotoreichik@ujf.cas.cz}
\maketitle
\begin{center}
	\vspace{-7ex}
	\textit{Dedicated to the memory of our friend and colleague Hagen Neidhardt}
\end{center}
\begin{abstract}
In this paper we deal with spectral optimization for the Robin Laplacian on a family of planar domains admitting parallel coordinates, namely a fixed-width strip built over a smooth closed curve and the exterior of a convex set with a smooth boundary. We show that if the curve length is kept fixed, the first eigenvalue referring to the fixed-width strip is for any value of the Robin parameter maximized by a circular annulus. Furthermore, we prove that the second eigenvalue in the exterior of a convex domain $\Omega$ corresponding to a negative Robin parameter does not exceed the analogous quantity for a disk whose boundary has a curvature larger than or equal to the maximum of that for $\partial\Omega$.
\end{abstract}
%


\section{Introduction}

\subsection{Motivation and the problem description}

Relations between geometry and spectral properties have a
rather long history. They are a trademark topic of mathematical
physics at least since the celebrated Faber and Krahn proof
\cite{F23, K24} of Lord Rayleigh's conjecture \cite{R77} about the
shape of the drum that produces the lowest tone. While the original
interest focused on problems with Dirichlet or Neumann boundary
conditions, more recently the attention shifted to Robin boundary
conditions. Spectral optimization for the Robin Laplacian was the
topic of a number of studies in the last few years and it still
offers many challenging open problems, see the reviews~\cite{BFK17,
L19} and the references therein.

It has to be said that spectral optimization can mean both the upper
and lower bounds. While in the mentioned Faber and Krahn result the
circle minimizes the principal eigenvalue, the Dirichlet Laplacian
on non-simply connected domains exhibits the opposite effect when
the full symmetry makes this eigenvalue maximal \cite{EHL99, HKK01}.
Moreover, this effect is robust, we note that a similar result holds
for a family of singular Schr\"odinger operators with an attractive
interaction supported by a closed planar curve \cite{EHL06}.

Recall that for a sufficiently regular planar domain $\Omega\subset\dR^2$ and the coupling constant $\alpha\in\dR$ the Robin eigenvalue problem can be written in PDE terms as
\[
	\begin{cases}
	-\Delta u = \lm u,&\quad\text{in}~\Omega,\\
	\frac{\partial u}{\partial \nu} + \alpha v = 0,&\quad\text{on}~\p\Omega
	\end{cases}
\]
where $\frac{\partial u}{\partial \nu}$ is the normal derivative of $u$ with the normal $\nu$ pointing outwards of $\Omega$. In the present paper we are going to deal with the optimization of the first and the second Robin eigenvalues on two particular classes of non-simply connected planar domains admitting the so-called parallel coordinates, specifically of loop-shaped curved strips and of exteriors of convex sets. These two classes of domains are closely related in the sense that the exteriors of a convex set can be viewed as the loop-shaped curved strip of infinite thickness built over the boundary. 

In the first case, we are going to prove that the lowest Robin eigenvalue on such a curved strip of a fixed length of the inner boundary and a fixed width is maximized by that of the annulus. This result can be regarded as an extension of the indicated property of the Dirichlet Laplacian~\cite[Thm~1a]{EHL99}. We stress that no restrictions are imposed here on the Robin coefficient, it can be negative as well as positive, in other words, the strip boundary can be both repulsive and attractive. The proof of the claim relies on the min-max principle, an appropriate test function being constructed via transplantation of the ground state for the annular strip using the parallel coordinates.

If the curve $\Sigma$ over which we build the domain is the boundary of a non-convex set, the existence of globally well defined parallel coordinates imposes a restriction on the strip width which could be expressed in terms of the \textit{absence of cut loci}, that is, points having the same distance from different parts of the domain boundary. Of course, any smooth curved strip of a fixed width has a cut locus represented by its axis, we have in mind  nontrivial ones referring to the distance from the curve $\Sigma$ only.

Our second result concerns an optimization of the second Robin eigenvalue in the exterior of a convex set under the assumption that the curvature $\kp$ of the boundary is non-negative and bounded from above by a fixed constant $\kp_\circ > 0$. The exterior of a disk with the boundary curvature $\kp_\circ$, in other words, with the radius $\kp_\circ^{-1}$, turns out be the unique maximizer. This new result complements the optimization of the lowest Robin eigenvalue in the exterior of a bounded set considered recently in~\cite{KL17a, KL17b}. It is not yet clear whether the stated condition on the curvature can be replaced by a more standard perimeter-type constraint. In contrast to the previous case, the Robin coefficient is now assumed to be negative as otherwise the spectrum of the respective Laplacian is purely essential and coincides with $[0,\infty)$, making thus the spectral optimization question void. In the proof, we again take advantage of the fact that the parallel coordinates are globally well defined. We apply the min-max principle to the span of the two transplanted eigenfunctions of the Robin Laplacian on the exterior of the disk corresponding to its first and second eigenvalues, respectively. Since the eigenfunction corresponding to the second eigenvalue is not radial, its transplantation is more involved and contains an additional geometric insight.

In the setting of bounded domains, it has been proved that the disk is a maximizer of the second Robin eigenvalue having a fixed area \cite{FL18a} (see also~\cite{FL18b}), provided that the negative boundary parameter lies in a specific interval. An analogous result has recently been proved in~\cite{GL19} for the third Robin eigenvalue with the maximizer being the union of two disjoint disks and with the negative boundary parameter again lying in a specific interval. In this context, we would like to emphasize that the optimization result for the second Robin eigenvalue in the present paper holds for \emph{all negative values} of the boundary parameter.

\subsection{Geometric setting}

Since the domain geometry is crucial in our results, let us first recall the necessary notions and state the assumptions we are going to use.

\begin{hypothesis}\label{hyp}
Let a $C^{\infty}$-smooth curve $\Sigma\subset\dR^2$ be the boundary of a bounded, simply connected domain $\Omega\subset\dR^2$. Let a circle $\cC\subset\dR^2$ be the boundary of a disk $\cB\subset\dR^2$. We denote by $L := |\Sigma|$ and $L_\circ := |\cC|$ the lengths of $\Sigma$ and $\cC$, respectively.
\end{hypothesis}
The mapping $\sigma\colon[0,L]\arr \dR^2$ provides the natural (counter-clockwise) parametrization of $\Sigma$ with the tangential vector $\tau(s) := \sigma'(s)$ satisfying $|\tau(s)| = 1$. We denote by $\kp \colon [0,L]\arr\dR$ the signed curvature of $\Sigma$; the convention we adopt is that $\kp\ge 0$ holds for convex $\Omega$.
Recall the \emph{Frenet formula}
\[
    \tau'(s) = -\kp\nu(s),
\]
where $\nu$ is the outer unit normal vector to the domain $\Omega$.

The object of our interest will be the curved strip built over $\Sigma$ with the thickness $d  \in (0,\infty]$, that is, the set
\begin{equation}\label{eq:strip}
    \Omega^{\rm c}_d
    :=
    \big\{x\in\dR^2\setminus\ov{\Omega}\colon
                    {\rm dist}\,(x,\Sigma) < d\big\}.
\end{equation}
The definition includes the unbounded domain $\Omega^{\rm c}_\infty$ identified with the exterior $\dR^2\sm\ov{\Omega}$ of $\Omega$ for which we will use the shorthand notation $\Omega^{\rm c} := \Omega^{\rm c}_\infty$. The boundary of $\Omega^{\rm c}_d$ is therefore
\[
\p\Omega^{\rm c}_d =
\begin{cases}
    \Sigma & \text{if}\;\: d = \infty,\\
    \Sigma\cup
    \big\{x\in\dR^2\setminus\ov{\Omega}\colon
    {\rm dist}\,(x,\Sigma) = d\big\} & \text{if}\;\: d < \infty.
\end{cases}
\]
In particular, $\p\Omega_d^{\rm c}$ has two components for $d < \infty$ and, respectively, one component for $d = \infty$. Consider the mapping
\begin{equation}\label{eq:mapping}
	[0,L)\tm(0,d)\ni (s,t)\mapsto \s(s)+t\nu(s).
\end{equation}
If the curvature $\kp \ge 0$ is sign-definite, then the mapping~\eqref{eq:mapping} is injective for all $d \in (0,\infty]$.
If the curvature $\kp$ is sign-changing, there exists by \cite[Prop. B.1]{BEHL17} a critical width $d_\star > 0$ such that~\eqref{eq:mapping} is injective for all $d < d_\star$ and in this case we assume that
\begin{equation}\label{eq:d}
    d \in (0, d_\star).
\end{equation}
The parallel coordinates $(s,t) \in [0,L)\times(0,d)$ on $\Omega^{\rm c}_d\,$ \cite{hart}, alternatively dubbed Fermi or natural curvilinear, are under the above assumption on $d$ everywhere well defined by the formula $\Omega^{\rm c}_d\ni x = \sigma(s) + t\nu(s)$.

\begin{remark}
As indicated in the introduction, the set of all $x\in\Omega^{\rm c}_d$ for which the closest point in $\p\Omega^{\rm c}_d$ is not uniquely defined is nonempty and coincides with $\{\s(s) + \frac{d}{2}\nu(s)\colon s\in[0,L)\}$. If the curvature $\kp$ is sign-changing, the existence of the parallel coordinates may be spoiled by a nontrivial cut locus referring to the distance from $\Sigma$ only, which may come from different sources, local and global. From \cite[Prop. B.1]{BEHL17}
we know that a necessary condition for the absence of such a cut locus is
\begin{equation}\label{eq:d}
    d\|\kappa_-\|_\infty < 1,
\end{equation}
where $\kappa_- = \min\{\kappa,0\}$. The condition \eqref{eq:d} is also sufficient, under the \textit{additional requirement} that $\int_s^{s'} \kp(s'')\,\dd s''>-\pi$ holds for all $0\le s<s'<L$. Indeed, the existence of a $\Sigma$-related cut locus would mean that the outer strip boundary must intersect itself, in other words, there must exist $s$ and $s'>s$ such that $\s(s)+d\nu(s) = \s(s')+d\nu(s')$. The curve $s\mapsto \s(\cdot)+d\nu(\cdot)$ is smooth, hence the angle between the tangents at the points with the parallel coordinates $(s,d)$ and $(s',d)$ cannot be larger than $-\pi$
(note that with our curvature convention the angles between tangents at nonconvex parts of the boundary are negative). 
However, the said tangents are parallel to the tangents $\tau(s)$ and $\tau(s')$ of the curve $\Sigma$, and since the angle between those is nothing else than $\int_s^{s'} \kp(s'')\,\dd s''$, the existence of such a cut locus leads to a contradiction, thus proving our claim.
\end{remark}
	
Since our setting is two-dimensional, it is useful to work with the complexified tangential and normal vectors
\begin{equation}\label{eq:t}
    {\bf t}(s) = \tau_1(s) + \ii\tau_2(s)
    \quad\;\text{and}\quad\;
    {\bf n}(s) = \nu_1(s) + \ii\nu_2(s).
\end{equation}
In this notation, the Frenet formula can be written in the complex form as
\begin{equation}\label{eq:Frenet}
    {\bf t}'(s) = -\kp{\bf n}(s)= \ii\kp{\bf t}(s).
\end{equation}

\subsection{The Robin Laplacian on $\Omega^{\rm c}_d$}

For an arbitrary value of the coefficient $\alpha\in\dR$, which characterizes the strength of the coupling to the boundary, we introduce the self-adjoint operator $\Op$ in the Hilbert space $L^2(\Omega^{\rm c}_d)$ through its quadratic form
\[
    \frm[u]
    :=
    \|\nabla u\|^2_{L^2(\Omega^{\rm c}_d;\dC^2)}
                + \alpha\|u|_{\p\Omega^{\rm c}_d}\|^2_{L^2(\p\Omega^{\rm c}_d)},
    \qquad \dom\frm = H^1(\Omega^{\rm c}_d),
\]
where $H^1(\Omega^{\rm c}_d)$ is the first-order $L^2$-based Sobolev space on $\Omega^{\rm c}_d$.

If $d < \infty$ the spectrum of $\Op$ is discrete and we denote by $\{\lm_k^\alpha(\Omega^{\rm c}_d)\}_{k\ge 1}$  its eigenvalues arranged in the non-decreasing order and repeated with multiplicities taken into account. The spectral properties of $\Opu$ corresponding to $d = \infty$ are different \cite[Prop. 2.1, Prop. 2.2]{KL17a}, namely
    \begin{myenum}
        \item $\sigma_{\rm ess}(\Opu) = [0,\infty)$.
        \item $\#\sigma_{\rm d}(\Opu) \ge 1$
        for all $\aa < 0$.
        \item $\sigma_{\rm d}(\Opu) = \varnothing$ for all $\aa \ge 0$.
        \end{myenum}
In analogy with the bounded domain case we denote by $\{\lm_k^\alpha(\Omega^{\rm c})\}_{k\ge 1}$ the negative eigenvalues of $\Opu$ arranged in the ascending order and repeated with the multiplicities taken into account. In the min-max spirit, this sequence is conventionally extended up to an infinite one by repeating the bottom of the essential spectrum $\inf\sigma_{\rm ess}(\Opu) = 0$ infinitely many times.

\subsection{Main results}

Let us now state our main results. The first one concerns optimization of $\lm_1^\aa(\Omega_d^{\rm c})$ on curved strips of a fixed width.
\begin{thm}\label{thm1}
    Assume that Hypothesis~\ref{hyp} holds and that $L = L_\circ$.
    Let $\kp$ be the curvature of $\Sigma$. The strip-width
    $d\in (0,\infty]$ may be arbitrary if $\kp \ge 0$, while for
    a sign-changing $\kp$ we assume that~\eqref{eq:d}
    is satisfied.
    Let the domains $\Omega_d^{\rm c}$ and $\cB^{\rm c}_d$ be as in~\eqref{eq:strip}. Then for
    the lowest Robin eigenvalues on these domains the inequality
    \[
        \lm_1^\alpha(\Omega_d^{\rm c})
        \le
        \lm_1^\alpha(\cB_d^{\rm c})
    \]
    holds for any $\alpha\in\dR$.
\end{thm}

Let us add a few comments. The above inequality holds trivially in the Neumann case, $\aa = 0$, since we have $\lm_1^0(\Omega_d^{\rm c}) = \lm_1^0(\cB_d^{\rm c}) = 0$. In the limit $\aa\arr+\infty$ it implies the respective inequality for the Dirichlet Laplacians providing thus an alternative proof of Theorem~1a in~\cite{EHL99}. Furthermore, if $\alpha < 0$, $\:\kp \ge 0$, and $d = \infty$, Theorem~\ref{thm1} reduces to the first claim of~\cite[Thm. 1.3]{KL17a}. Note that the geometric character of $\Omega$ manifested in the constraint on the `thickness' plays a role again: the annulus is always a maximizer here, even for $\alpha > 0$, while in the case of general bounded domains under fixed area constraint the disk is conjectured to be a maximizer in the subclass of simply-connected domains for $\aa < 0$ and is known to be a minimizer for $\aa > 0$, \cf~\cite{FK15, AFK17} in the former case and \cite{Bossel_1986, Daners_2006} in the latter. Moreover, for general bounded domains the disk is a maximizer for $\aa < 0$ under fixed perimeter constraint~\cite{AFK17}. Multi-dimensional analogues of the latter result are obtained in \cite{BFNT18, V19}.

As already mentioned in the introduction, the proof of Theorem~\ref{thm1} will rely on the min-max principle with a suitable test function constructed through the transplantation of the radial ground-state eigenfunction for the annulus using the method of parallel coordinates.

Our second result concerns optimization of $\lm_2^\aa(\Omega^{\rm c})$ on unbounded exterior domains described above.
\begin{thm}\label{thm2}
    Assume that Hypothesis~\ref{hyp} holds.
    Let $\kp\colon[0,L]\arr\dR$ and $\kp_\circ \in\dR_+$ be the curvatures of $\Sigma$ and $\cC$, respectively.
    Assume further that $\Omega$ is convex, that is, $\kp \ge 0$, and that $\max\kp \le \kp_\circ$ holds .
    Let the domains $\Omega^{\rm c}$ and $\cB^{\rm c}$ be as in~\eqref{eq:strip}
    with $d = \infty$. Then for the second Robin eigenvalues on these domains the inequality
    \begin{equation}\label{eq:main2}
        \lm_2^\alpha(\Omega^{\rm c})
        \le
        \lm_2^\alpha(\cB^{\rm c})
    \end{equation}
    is valid for any $\alpha < 0$. If $\lm_2^\alpha(\cB^{\rm c}) < 0$ and
    the equality in~\eqref{eq:main2} holds, the two domains are congruent, $\Omega\cong\cB$.
\end{thm}
The above theorem and monotonicity of
$\lm_2^\alpha(\cB^{\rm c})$, with respect to $L_\circ$ shown in Proposition~\ref{prop:annulus2}, below yield the
following.
\begin{cor}\label{cor}
    Assume that Hypothesis~\ref{hyp} holds and let $\kp_\circ > 0$ be fixed.
    Then, for all $\aa < 0$,
    \begin{equation}\label{eq:main3}
    \max_{\stackrel{\Omega~{\rm convex}}{\kp \le \kp_\circ}}\lm_2^\alpha(\Omega^{\rm c})
    =
    \lm_2^\alpha(\cB^{\rm c}),
    \end{equation}
    where the maximum is taken over all
    convex smooth domains $\Omega\subset\dR^2$
    whose curvature satisfies $\max\kp \le \kp_\circ$ and where $\cB\subset\dR^2$ is a disk of the curvature $\kp_\circ$.
\end{cor}

We remark that the inequality~\eqref{eq:main2} is nontrivial only if $\sfH_{\aa,\cB^{\rm c}}$ has more than one negative eigenvalue. We also emphasize that, in contrast to Theorem~\ref{thm1} we have $L \ne L_\circ$ in general, in fact, it is easy to show that $L > L_\circ$ holds unless $\Omega \cong\cB$. In order to prove Theorem~\ref{thm2}, we apply the min-max principle transplanting to $\Omega^{\rm c}$ the span of the two eigenfunctions of $\sfH_{\aa,\cB^{\rm c}}$ corresponding to the eigenvalues $\lm_1^\aa(\cB^{\rm c})$ and $\lm_2^\aa(\cB^{\rm c})$, respectively. The ground-state is transplanted in a conventional way, however, the transplantation of the first excited state is a little more involved.
We note that an eigenfunction corresponding to the second Robin eigenvalue on the exterior of a disk can be written in parallel coordinates on $\cB^{\rm c}$ as
\[
    v_\circ(s,t) =
    \phi(t)\exp\left(\frac{2\pi\ii}{L_\circ}s\right).
\]
Since $\exp\left(\frac{2\pi\ii}{L_\circ}s\right)$ can be interpreted as the complexified tangent vector for $\cB$, a natural way of transplantation of $v_\circ$ onto $\Omega^{\rm c}$ would be
\[
    v_\star(s,t) = \phi(t) {\bf t}(s),
\]
where ${\bf t}$ is the complexified tangent vector for $\Omega$ defined in~\eqref{eq:t}.

\section{Preliminaries}

\subsection{The quadratic form $\frm$ in parallel coordinates}

Our first main tool is the representation of the quadratic form $\frm$ in the parallel coordinates on $\Omega^{\rm c}_d$. Using them, the inner product in the Hilbert space $L^2(\Omega_d^{\rm c})$ can be written as follows,
\[
    (u,v)_{L^2(\Omega_d^{\rm c})}
    =\int_0^d\int_0^L u(s,t)\ov{v(s,t)}
    \big(1+\kp(s) t\big)\,\dd s\,\dd t.
\]
It is well known that the gradient in these coordinates is expressed as
\[
    \nabla u =
    \frac{\tau(s)}{1+\kp(s) t}\,\p_s u +
    \nu(s)\p_t u.
\]
Consequently, the quadratic form $\frm$ can be written in the parallel coordinates for $d < \infty$ as
\begin{equation}\label{eq:form_parallel}
\begin{aligned}
    \frm[u]
    &\! = \!\int_0^d\int_0^L\left(
    \frac{|\p_s u(s,t)|^2}{1+\kp(s)t}+
    |\p_t u(s,t)|^2(1+\kp(s)t)
    \right)\,\dd s\,\dd t\\
    &\qquad\qquad + \aa\int_0^L |u(s,0)|^2\,\dd s + \aa\int_0^L |u(s,d)|^2(1+\kp(s)d)\,\dd s,\\
    \dom\frm &\! =\!
    \left\{u\colon \Sigma\tm (0,d)\arr\dC\colon\!
\int_0^d\!\int_0^L\left[
\frac{|\p_s u|^2}{1+\kp t}+
(|u|^2\!+\! |\p_t u|^2)(1+\kp t)
\right]\,\dd s\,\dd t < \infty\right\}.
\end{aligned}
\end{equation}
The above representation remains valid for $d = \infty$ without the last term corresponding to the outer boundary, provided that $\Omega$ is convex.

\subsection{Eigenfunctions in the radially symmetric case}

We also need properties of the eigenfunctions corresponding to the first and the second eigenvalue in the radially symmetric case. They are elementary but we describe them in the next two propositions, the proofs of which are postponed to the appendices, in order to make the paper self-contained. Let us begin with the ground-state eigenfunction of the Robin annulus.
\begin{prop}\label{prop:annulus1}
    Assume that Hypothesis~\ref{hyp} holds.
    For any fixed $d > 0$ and any $\alpha\in\dR$,
    or for $d = \infty$ and any $\aa <0$,
    the lowest  eigenvalue $\lm_1^\aa(\cB^{\rm c}_d)$ of
    $\OpD$ is simple and
    the corresponding eigenfunction
    can be written
    in the parallel coordinates on $\cB_d^{\rm c}$
    as
    \[
    u_\circ(s,t) = \psi(t),
    \]
    with a given real-valued
    $\psi \in C^\infty([0,d])$ if $d < \infty$
    and with $\psi \in C^\infty([0,\infty))$
    satisfying
    \begin{equation}\label{eq:integrability1}
    \int_0^\infty\big[\psi(t)^2 + \psi'(t)^2](1+t)\,\dd t < \infty,
    \end{equation}
    if $d = \infty$.
\end{prop}

Consider next the first excited state of the Robin Laplacian in the exterior of a disk.
\begin{prop}\label{prop:annulus2}
    Assume that Hypothesis~\ref{hyp} holds.
    Then for any fixed $\alpha < 0$ such that
    $\#\s_{\rm d}(\OpDu) > 1$,
    the second eigenvalue $\lm_2^\aa(\cB^{\rm c}) <0$ of $\OpDu$
    has multiplicity two
    and the respective eigenfunctions of
    $\OpDu$ can be written
    in parallel coordinates on $\cB^{\rm c}$
    as
    \[
    v_\circ^\pm(s,t) = \exp\left(
    \pm \frac{2\pi\ii }{L_\circ}s\right)
    \phi(t),\quad\; s\in [0,L_\circ),\; t\in [0,\infty),
    \]
    with a given real-valued
    $\phi\in C^\infty([0,\infty))$
    satisfying the integrability condition
    \begin{equation}\label{eq:integrability2}
    \int_0^\infty\big[\phi(t)^2 + \phi'(t)^2](1+t)\,\dd t < \infty.
    \end{equation}
    Moreover, $\lm_2^\alpha(\cB^{\rm c})$
        is a non-increasing function
    of $L_\circ$.
\end{prop}

We remark that the functions $\psi$ and $\phi$ in Propositions~\ref{prop:annulus1} and~\ref{prop:annulus2} can be explicitly expressed in terms of Bessel functions, however, this is not essential for our analysis.

\section{Proofs of the main results}

Now we are going to provide proofs of Theorems~\ref{thm1} and~\ref{thm2}. Recall that the $C^\infty$-smooth curve $\Sigma\subset\dR^2$ is the boundary of a bounded, simply connected domain $\Omega\subset\dR^2$, and the circle $\cC\subset\dR^2$ is the boundary of the disk $\cB\subset\dR^2$. The lengths of $\Sigma$ and $\cC$ are denoted by $L$ and $L_\circ$, respectively. The curvature of $\Sigma$ is denoted by $\kp$ and the curvature of $\cC$ is a constant $\kp_\circ > 0$.

\subsection{Proof of Theorem~\ref{thm1}}

By assumption we have $L = L_\circ$ and we fix $d > 0$ satisfying the additional condition \eqref{eq:d} in the case that $\kp$ is sign-changing. Furthermore, $\alpha\in\dR$ is an arbitrary fixed number. The case $d = \infty$ is dealt with in~\cite[Thm~1.3]{KL17a} and thus we may omit it here. By Proposition~\ref{prop:annulus1}, there exists a function $\psi \in C^\infty([0,d])$ such that the ground-state $u_\circ \in C^\infty(\Omega^{\rm c}_d)$ of $\OpD$ can be written as $u_\circ(s,t) = \psi(t)$ in the parallel coordinates on $\cB^{\rm c}_d$. Using it we define the test function $u_\star\in H^1(\Omega_d^{\rm c})$ in the parallel coordinates on the curved strip $\Omega^{\rm c}_d$ as follows,
\[
    u_\star(s,t) := \psi(t),\qquad
    s\in [0,L],\,   t\in [0,d].
\]
Using the representation of $\frm$ in~\eqref{eq:form_parallel}, applying the min-max principle and the total curvature identity $\int_0^L \kp(s)\dd s = 2\pi$  we obtain
\[
\begin{aligned}
    \lm_1^\alpha(\Omega^{\rm c}_d)
    & \le \frac{\frm[u_\star]}{\|u_\star\|^2_{L^2(\Omega^{\rm c}_d)}}\\
    & =
    \frac{\displaystyle
        \int_0^d\int_0^L\psi'(t)^2(1+ \kp(s) t)\,\dd s\,\dd t
        +
        \alpha \int_0^L\Big[|\psi(0)|^2 + |\psi(d)|^2(1+d\kp(s))\Big] \,\dd s}{\displaystyle    \int_0^d\int_0^L\psi(t)^2(1+ \kp(s) t)\,\dd s\,\dd t}
    \\
    & =
    \frac{\displaystyle
    \int_0^d\psi'(t)^2(L+ 2\pi t)\,\dd t
    +
    \alpha L |\psi(0)|^2 + \alpha (L+2\pi d)|\psi(d)|^2}{\displaystyle  \int_0^d\psi(t)^2(L+ 2\pi t)\,\dd t}\\
    & =
    \frac{\frmD[u_\circ]}{\|u_\circ\|^2_{L^2(\cB^{\rm c}_d)}}
    = \lm_1^\alpha(\cB^{\rm c}_d),
\end{aligned}
\]
which yields the sought claim.

\subsection{Proof of Theorem~\ref{thm2}}

In view of the convexity of $\Omega$ the curvature of $\Sigma$ satisfies $\kp\ge 0$ and by assumption $\max\kp \le  \kp_\circ$ holds. Let us exclude the trivial case supposing that  $\Omega\ncong\cB$. Then we have $\min\kp < \kp_\circ$ which implies
\begin{equation}\label{eq:L}
    L = \frac{L\kp_\circ}{\kp_\circ} > \frac{\displaystyle\int_0^L\kp(s)\,\dd s}{\kp_\circ}
    = \frac{2\pi}{\kp_\circ} = L_\circ.
\end{equation}
We fix the `width' $d =\infty$ and the coupling constant $\alpha < 0$. Without loss of generality, we may assume that $|\aa|$ is large enough so that $\lm_2^\aa(\cB^{\rm c}) < 0$
as otherwise the inequality~\eqref{eq:main2} would trivially hold.

\medskip

\noindent {\bf Step 1.} \emph{Test functions}.
In view of Propositions~\ref{prop:annulus1} and~\ref{prop:annulus2}, we can represent the eigenfunctions of $\OpDu$ corresponding to its simple first eigenvalue $\lm_1^\alpha(\cB^{\rm c})$ and the second eigenvalue $\lm_2^\alpha(\cB^{\rm c})$ of multiplicity two in parallel coordinates $(s,t)$ on $\cB^{\rm c}$ as
\begin{equation}\label{eq:EFs}
    u_\circ(s,t) = \psi(t)
    \qquad\text{and}\qquad
    v_\circ^\pm(s,t) =
    \exp\left(\pm\frac{2\pi\ii s}{L_\circ}
    \right)\phi(t),
\end{equation}
where $\psi,\phi\in C^\infty([0,\infty))$ are real-valued and satisfy the integrability conditions~\eqref{eq:integrability1} and~\eqref{eq:integrability2}, respectively. We introduce test functions $u_\star,v_\star \in H^1(\Omega^{\rm c})$ on $\Omega^{\rm c}$ defining  them in terms of the parallel coordinates as
\[
    u_\star(s,t) := \psi(t)\qquad\text{and}\qquad
    v_\star(s,t) := {\bf t}(s)\phi(t),
    \qquad s\in[0,L],\, t\in [0,\infty),
\]
where ${\bf t}(s)$ is the complexified normal~\eqref{eq:t}.

\medskip

\noindent {\bf Step 2.} \emph{Orthogonality.}
Next, we are going to show that $u_\star$ and $v_\star$ are orthogonal in $L^2(\Omega^{\rm c})$.
 To this aim, we observe that
\[
    \int_0^L {\bf t}(s)\,\dd s
    =
    \int_0^L (\sigma_1'(s) + \ii \sigma_2'(s))\,\dd s =
    \sigma_1(L) + \ii\sigma_2(L) -
    \sigma_1(0) - \ii\sigma_2(0)  =0,
\]
%
where the fact that $\Sigma$ is a closed curve was employed.
Furthermore, using the Frenet formula we get
\[
    \int_0^L {\bf t}(s)\kp(s)\,\dd s
    =
    -\ii\int_0^L {\bf t}'(s)\,\dd s
    = -\ii({\bf t}(L) - {\bf t}(0)) = 0,
\]
where the closedness and smoothness of $\Sigma$ were taken into account.
Combining these two relations we infer that
\begin{equation}\label{eq:ortho1}
\begin{aligned}
    (v_\star,u_\star)_{L^2(\Omega^{\rm c})}
    & \!= \!
    \int_0^\infty\int_0^L\psi(t)\phi(t)
    {\bf t}(s)(1+ t\kp(s))\,\dd s \,\dd t\\
    & \! = \!
    \int_0^\infty\int_0^L\psi(t)\phi(t){\bf t}(s)\,\dd s \,\dd t\! +\!
    \int_0^\infty\int_0^Lt\psi(t)\phi(t)
    {\bf t}(s)\kp(s)\,\dd s \,\dd t
     \!=\! 0.
\end{aligned}
\end{equation}
At the same time, we have
\begin{equation}\label{eq:ortho2}
\begin{aligned}
    \frmu[v_\star,\!u_\star] \! = \!
    \int_0^\infty\!\int_0^L\psi'(t)\phi'(t)
    {\bf t}(s)(1\!+\! t\kp(s))\,\dd s \,\dd t
     \!+\!
    \alpha\psi(0)\phi(0)\!
    \int_0^L
    {\bf t}(s)\,\dd s\! =\! 0.
\end{aligned}
\end{equation}

\noindent {\bf Step 3.} \emph{Bounds on the Rayleigh
    quotients.}
For a non-trivial function $u \in H^1(\Omega^{\rm c})$ we define
\[
    \sfR_{\aa,\Omega^{\rm c}}[u] :=
        \frac{\frmu[u]}{\|u
        \|^2_{L^2(\Omega^{\rm c})}}.
\]
Using~\eqref{eq:form_parallel} and the total curvature identity $\int_0^L \kp(s)\,\dd s = 2\pi$, the Rayleigh quotient of the test function $u_\star$ defined in this way
can be expressed as
\[
\begin{aligned}
    \sfR_{\aa,\Omega^{\rm c}}[u_\star]
    & =
    \frac{\displaystyle\int_0^\infty\int_0^L\psi'(t)^2(1+t\kp(s))
        \,\dd s
        \,\dd t + \aa L|\psi(0)|^2}{\displaystyle\int_0^\infty\int_0^L\psi(t)^2(1+t\kp(s))
        \,\dd s\,\dd t}\\
    &  =
    \frac{\displaystyle\int_0^\infty
        \psi'(t)^2(L+2\pi t)
        \,\dd t + \aa L|\psi(0)|^2}{\displaystyle\int_0^\infty\psi(t)^2(L+2\pi t)\,\dd t}\\
    & =
        \frac{\displaystyle\int_0^\infty\psi'(t)^2\left(1+\frac{2\pi t}{L}\right)
        \,\dd t + \aa|\psi(0)|^2}{\displaystyle\int_0^\infty\psi(t)^2\left(1+\frac{2\pi t}{L}\right)\,\dd t}.
\end{aligned}
\]
Furthermore, using the inequalities $\lm_1^\aa(\cB^{\rm c}) < 0$ and $L > L_\circ$, we get the following estimate
\begin{equation}\label{eq:bound1}
    \sfR_{\aa,\Omega^{\rm c}}[u_\star]
     \le
    \frac{\displaystyle\int_0^\infty\psi'(t)^2\left(1+\frac{2\pi t}{L_\circ}\right)
    \,\dd t + \aa|\psi(0)|^2}{\displaystyle\int_0^\infty\psi(t)^2\left(1+\frac{2\pi t}{L_\circ}\right)\,\dd t}
    =\mathsf{R}_{\alpha,\cB^{\rm c}}[u_\circ]
    =\lm_1^\aa(\cB^{\rm c}).
\end{equation}
Making use of~\eqref{eq:form_parallel}, the total curvature identity and the Frenet formula~\eqref{eq:Frenet}, the Rayleigh quotient corresponding to
the test function $v_\star$ is given by
\[
\begin{aligned}
    \mathsf{R}_{\aa,\Omega^{\rm c}}[v_\star]
    &\! =\!
    \frac{\displaystyle\int_0^\infty\!\int_0^L\phi'(t)^2(1\!+\!t\kp(s))
    \,\dd s
    \,\dd t
    \!+\!
    \int_0^\infty\!\int_0^L
    \frac{\kp^2(s)\phi(t)^2}{1+t\kp(s)}
    \,\dd s
    \,\dd t
    \!+\!
    \aa L|\phi(0)|^2}{\displaystyle\int_0^\infty\int_0^L\phi(t)^2(1+t\kp(s))
    \,\dd s\,\dd t}\\
    &\! =\!
    \frac{\displaystyle\int_0^\infty\phi'(t)^2(L+2\pi t)\,\dd t
        \!+\!
        \int_0^\infty\!\int_0^L
        \frac{\kp^2(s)\phi(t)^2}{1+t\kp(s)}
        \,\dd s
        \,\dd t
        \!+\!
        \aa L|\phi(0)|^2}{\displaystyle\int_0^\infty
        \phi(t)^2(L+2\pi t)\,\dd t}.
\end{aligned}
\]
Using further the strict monotonicity of the function
\[
    \dR_+\ni x\mapsto \frac{x^2}{1+ tx},\quad\; t \ge0,
\]
in combination with the inequalities $L > L_\circ$, $\,\max\kp\le \kp_\circ$, and $\min\kp < \kp_\circ$,
we get for the Rayleigh quotient corresponding to $v_\star$ the following estimate,
\begin{equation}\label{eq:bound2}
\begin{aligned}
\mathsf{R}_{\aa,\Omega^{\rm c}}[v_\star]
& <
\frac{\displaystyle\int_0^\infty\phi'(t)^2(L+2\pi t)
    \,\dd t +
    L\int_0^\infty
    \frac{\kp^2_\circ   \phi(t)^2}{1+t\kp_\circ}
    \,\dd t
    +
    \aa L |\phi(0)|^2}{\displaystyle\int_0^\infty\phi(t)^2(L+
    2\pi t)\,\dd t}\\
& =
    \frac{\displaystyle\int_0^\infty\phi'(t)^2\left(1+\frac{2\pi t}{L}\right)
    \,\dd t
    +
    \int_0^\infty
    \frac{\kp^2_\circ   \phi(t)^2}{1+t\kp_\circ}
    \,\dd t
    +
    \aa  |\phi(0)|^2}{\displaystyle\int_0^\infty\phi(t)^2\left(1+
    \frac{2\pi t}{L}\right)\,\dd t}\\
    & \le
    \frac{\displaystyle\int_0^\infty\phi'(t)^2\left(1+\frac{2\pi t}{L_\circ}\right)
    \,\dd t
    +
    \int_0^\infty
    \frac{\kp^2_\circ\phi(t)^2}{1+t\kp_\circ}
    \,\dd t
        +
        \aa  |\phi(0)|^2}{\displaystyle\int_0^\infty\phi(t)^2\left(1+
        \frac{2\pi t}{L_\circ}\right)\,\dd t}\\
        & =
        \sfR_{\alpha,\cB^{\rm c}}[v_\circ]
        =   \lm_2^\aa(\cB^{\rm c}).
\end{aligned}
\end{equation}

\noindent {\bf Step 4.} \emph{The min-max principle.} Any $w_\star\in{\rm span}\,\{u_\star,v_\star\}\sm\{0\}$ can be represented as a linear combination $w_\star = p u_\star + qv_\star$ with $(p,q)\in\dC^2_\tm:=\dC^2\sm\{(0,0)\}$. The following simple inequality,
\begin{equation}\label{eq:trivial_ineq}
    \frac{a+b}{c+d} \le \max\left\{\frac{a}{c},\frac{b}{d}\right\},
\end{equation}
holds obviously for any $a,b \in\dR$ and $c,d > 0$. Applying the min-max principle, using the orthogonality relatons~\eqref{eq:ortho1},~\eqref{eq:ortho2}, the bounds~\eqref{eq:bound1},~\eqref{eq:bound2}, and the inequality~\eqref{eq:trivial_ineq} we get
\[
\begin{aligned}
    \lm_2^\alpha(\Omega^{\rm c})
    &\le
    \max_{(p,q)\in\dC^2_\tm}
    \frac{\frmu[p u_\star + qv_\star]}{\|p u_\star + qv_\star\|^2_{L^2(\Omega^{\rm c})}}\\
    &
    \hspace{-2em} = \max_{(p,q)\in\dC^2_\tm}
    \frac{|p|^2\frmu[ u_\star]
    +|q|^2\frmu[v_\star]}{
    |p|^2\|u_\star\|^2_{L^2(\Omega^{\rm c})}
    + |q|^2\|v_\star\|^2_{L^2(\Omega^{\rm c})}}
    \le
    \max\left\{
    \sfR_{\aa,\Omega^{\rm c}}[u_\star],
    \sfR_{\aa,\Omega^{\rm c}}[v_\star]
    \right\}
    <
    \lm_2^\alpha(\cB^{\rm c}),
\end{aligned}
\]

\subsection*{Acknowledgment}

The research was supported by the Czech Science Foundation within the project 17-01706S. P.E. also acknowledges support of the EU project \\ CZ.02.1.01/0.0/0.0/16\textunderscore 019/0000778. The authors are indebted to Magda Khalile for fruitful discussions and to the referees for useful comments.

\begin{appendix}

\section{Proof of Proposition~\ref{prop:annulus1}}

The case $d = \infty$ was dealt with in~\cite[Sec. 3]{KL17a}. Assume that $d < \infty$ and let $\aa\in\dR$ be arbitrary. For the sake of simplicity and without loss of generality, we also assume that $L_\circ = 2\pi$. In this case the curvilinear coordinates $(s,t)$ essentially coincide with the polar coordinates. Using the complete family of orthogonal projections on $L^2(\cB_d^{\rm c})$,
\[
(\Pi_n u)(s,t) = \frac{1}{\sqrt{2\pi}}
\ee^{\ii n s}\int_0^{2\pi}
u(s',t)\,\frac{\ee^{-\ii n s'}}{\sqrt{2\pi}}\,\dd s',\quad\; n\in\dZ,
\]
one can decompose $\OpD$ into an orthogonal sum
\[
\OpD = \bigoplus_{n\in\dZ}\OpD^{[n]},
\]
where the self-adjoint fiber operator $\OpD^{[n]}$ acts in the Hilbert space $L^2((0,d);(1+ t)\,\dd t)$ and corresponds to the quadratic form $(n\in\dZ)$
\[
\begin{aligned}
\frmD^{[n]}[\psi]
& =\int_0^d
\left(|\psi'(t)|^2(1+t) +\frac{n^2|\psi(t)|^2}{1+t}\right)\,\dd t + \aa|\psi(0)|^2 +\aa(1+d)|\psi(d)|^2,\\
\dom\frmD^{[n]} &= H^1((0,d)).
\end{aligned}
\]
Clearly the lowest eigenvalue of $\OpD^{[0]}$ is simple and strictly smaller than the lowest eigenvalues of $\OpD^{[n]}$ with $n\ne0$. Thus, the ground-state of $\OpD$ is simple and depends on $t$ variable only. The smoothness of the corresponding eigenfunction follows from standard elliptic regularity theory.

\section{Proof of Proposition~\ref{prop:annulus2}}

Using the complete family of orthogonal projections on $L^2(\cB^{\rm c})$
\[
(\Pi_n u)(s,t) = \frac{1}{\sqrt{L_\circ}}
\ee^{\frac{2\pi\ii n s}{L_\circ}}\int_0^{L_\circ}
u(s',t)\,\frac{\ee^{-\frac{2\pi\ii n s'}{L_\circ}}}{\sqrt{L_\circ}}\,\dd s',\quad\; n\in\dZ,
\]
one can again decompose $\OpDu$ into an orthogonal sum
\[
\OpDu = \bigoplus_{n\in\dZ} \OpDu^{[n]},
\]
where the fiber operators $\OpDu^{[n]}$, $\:n\in\dZ$, in the Hilbert space $L^2(\dR_+;(1 + \frac{2\pi t}{L_\circ})\,\dd t)$ correspond to the quadratic forms
\[
\begin{aligned}
\frmDu^{[n]}[\psi]
& =\int_0^\infty
\left(|\psi'(t)|^2\left(1 + \frac{2\pi t}{L_\circ}\right) +\frac{1}{L_\circ}\frac{4\pi^2n^2|\psi(t)|^2}{L_\circ + 2\pi t}\right)\,\dd t + \aa |\psi(0)|^2,\\
\dom\frmDu^{[n]} &=
\left\{\psi\colon\dR_+\arr\dC\colon
\psi,\psi'\in L^2(\dR_+;\left(1 + 2\pi L_\circ^{-1} t\right)\,\dd t)\right\}.
\end{aligned}
\]
It is easy to see that $\OpDu^{[0]}$ has exactly one negative simple eigenvalue, which corresponds to the ground-state eigenvalue $\lm_1^\aa(\cB^{\rm c})$ of $\OpDu$. The first excited state eigenvalue $\lm_2^\aa(\cB^{\rm c})$ corresponds to the lowest eigenvalues of the identical operators $\OpDu^{[1]}$ and $\OpDu^{[-1]}$. Moreover, the smoothness of $\phi$ follows from standard elliptic regularity theory.

Let $\cB_1,\cB_2$ be two disks with
 perimeters
$L_1$ and $L_2$, respectively. Assume that
$L_1< L_2$. Then we obtain that
\[
    \begin{aligned}
    \lm_2^\aa(\cB^{\rm c}_1)
    & =
    \inf_{\psi\in C^\infty_0([0,\infty))}
    \frac{\displaystyle
        \int_0^\infty
        \left(|\psi'(t)|^2\left(1 + \frac{2\pi}{L_1} t\right) +\frac{1}{L_1}
        \frac{4\pi^2|\psi(t)|^2}{L_1 + 2\pi t}\right)\,\dd t + \aa |\psi(0)|^2}{
        \displaystyle
        \int_0^\infty
    |\psi(t)|^2\left(1 + \frac{2\pi}{L_1} t\right)\dd t} \\
    & \ge
    \inf_{\psi\in C^\infty_0([0,\infty))}
    \frac{\displaystyle
        \int_0^\infty
        \left(|\psi'(t)|^2\left(1 + \frac{2\pi}{L_2} t\right) +\frac{1}{L_2}
        \frac{4\pi^2|\psi(t)|^2}{L_2 + 2\pi t}\right)\,\dd t + \aa |\psi(0)|^2}{
        \displaystyle
        \int_0^\infty
        |\psi(t)|^2\left(1 + \frac{2\pi}{L_2} t\right)\dd t} \\
    & =
        \lm_2^\aa(\cB^{\rm c}_2).
    \end{aligned}
\]
Hence, it follows that $\lm_2^\alpha(\cB^{\rm c})$ is a non-increasing function of its perimeter.

\end{appendix}

%

%

\newcommand{\etalchar}[1]{$^{#1}$}

\end{document}

%% file: commands.tex
\usepackage{scrextend}
\changefontsizes{10.5pt}

\newcommand\frm{\frh_{\alpha,\Omega^{\rm c}_d}}
\newcommand\frmu{\frh_{\alpha,\Omega^{\rm c}}}
\newcommand\frmD{\frh_{\alpha,\cB^{\rm c}_d}}
\newcommand\frmDu{\frh_{\alpha,\cB^{\rm c}}}

\newcommand\Op{\sfH_{\alpha,\Omega^{\rm c}_d}}
\newcommand\OpD{\sfH_{\alpha,\cB^{\rm c}_d}}
\newcommand\Opu{\sfH_{\alpha,\Omega^{\rm c}}}
\newcommand\OpDu{\sfH_{\alpha,\cB^{\rm c}}}

\newcommand\s{\sigma}

\newcommand\frh{\mathfrak{h}}

\newcommand{\cf}{{\it cf.}\,}

\renewcommand\and{\qquad\text{and}\qquad}

\newcommand\sm{\setminus}
\newcommand\tm{\times}

\newcommand{\comm}[1]{}

\renewcommand\aa{\alpha}
\newcommand\lm{\lambda}

\newcommand\p{\partial}

\newcommand\kp{\kappa}

\newcommand\ii{{\mathsf{i}}}

\newcommand\arr{\rightarrow}

\newcommand{\ee}{\mathsf{e}}

\newcommand\dd{{\mathsf{d}}}

\newcounter{counter_a}
\newenvironment{myenum}{\begin{list}{{\rm(\roman{counter_a})}}%
{\usecounter{counter_a}
\setlength{\itemsep}{1.ex}\setlength{\topsep}{0.8ex}
\setlength{\leftmargin}{5ex}\setlength{\labelwidth}{5ex}}}{\end{list}}

\usepackage[latin1]{inputenc}
\usepackage[T1]{fontenc}

\numberwithin{figure}{section}
\numberwithin{equation}{section}
\theoremstyle{plain}
\newtheorem*{thm*}{Theorem}
\newtheorem{thm}{Theorem}[section]

\newtheorem{prop}[thm]{Proposition}

\newtheorem{cor}[thm]{Corollary}

\theoremstyle{remark}
\newtheorem{remark}[thm]{Remark}
\theoremstyle{plain}
\newtheorem{hypothesis}{Hypothesis}

\newcommand\ov{\overline}

\def\ov{\overline}

      \def\dC{{\mathbb C}}

      \def\dR{{\mathbb R}}

   \def\dZ{{\mathbb Z}}

   \def\sfH{{\mathsf H}}

   \def\cB{{\mathcal B}}   \def\cC{{\mathcal C}}

      \def\sfR{{\mathsf R}}

\newcommand{\dom}{\mathrm{dom}\,}

\makeatletter
\def\section{\@startsection{section}{1}\z@{.9\linespacing\@plus\linespacing}%
	{.7\linespacing} {\fontsize{13}{14}\selectfont\bfseries\centering}}
\def\paragraph{\@startsection{paragraph}{4}%
	\z@{0.3em}{-.5em}%
	{$\bullet$ \ \normalfont\itshape}}

\@addtoreset{equation}{section}
\makeatother

%% file: v_final.bbl
\begin{thebibliography}{\textsc{BCD{\etalchar{+}}72}}


\bibitem[AFK17]{AFK17}
P.~R.~S. Antunes, P.~Freitas, and D.~Krej\v{c}i\v{r}\'{\i}k, Bounds and extremal domains for Robin eigenvalues with negative boundary parameter,
{\it Adv. Calc. Var.} \textbf{10} (2017), 357--380.

\bibitem[BEHL17]{BEHL17}
J.~Behrndt, P.~Exner, M.~Holzmann, and V.~Lotoreichik,
Approximation of Schr\"{o}dinger operators with $\delta$-interactions supported on hypersurfaces,
{\it Math. Nachr.} {\bf 290} (2017), 1215--1248.

\bibitem[Bos86]{Bossel_1986}
M.-H. Bossel,
Membranes {\'{e}}lastiquement li{\'{e}}es: {E}xtension du th{\'{e}}or{\'{e}}me de {R}ayleigh-{F}aber-{K}rahn et de l'in{\'{e}}galit{\'{e}} de {C}heeger,
{\it C. R. Acad. Sci. Paris S\'{e}r. Math.} \textbf{302} (1986), 47--50.

\bibitem[BFNT19]{BFNT18}
D. Bucur, V. Ferone, C. Nitsch, and  C. Trombetti, A sharp estimate for the first Robin-Laplacian eigenvalue with negative boundary parameter, {\it Atti Accad. Naz. Lincei, Cl. Sci. Fis. Mat. Nat., IX. Ser., Rend. Lincei, Mat. Appl.,} {\bf 30} (2019), 665--676.

\bibitem[BFK17]{BFK17}
D. Bucur, P. Freitas, and J. B. Kennedy: The Robin problem, in `Shape optimization and spectral theory' (A. Henrot, ed.), De Gruyter 2017; pp. 78-119.

%
%
%
%
%



\bibitem[Dan06]{Daners_2006}
D.~Daners, A Faber-Krahn inequality for Robin problems in any space dimension, {\it Math. Ann.} \textbf{335} (2006), 767--785.


\bibitem[EHL99]{EHL99}
P. Exner, E.\,M. Harrell, and M. Loss,
Optimal eigenvalues for some Laplacians and Schroedinger operators depending on curvature, {\it Oper. Theory: Adv. Appl.} \textbf{108} (1999), 47--58.

\bibitem[EHL06]{EHL06}
P.~Exner, E.\,M.~Harrell, and M.~Loss, Inequalities for means of chords, with application to isoperimetric problems,
{\it Lett. Math. Phys.} \textbf{75} (2006), 225-233; addendum \textbf{77} (2006), 219.

\bibitem[F23]{F23}
G. Faber, Beweis  dass  unter  allen  homogenen  Membranen  von  gleicher  Fl\"ache  und gleicher  Spannung  die  kreisf\"ormige  den  tiefsten  Grundton  gibt, {\it Sitz.  bayer.  Akad.Wiss.} (1923), 169--172.
%


\bibitem[FK15]{FK15}
P.~Freitas and D.~Krej\v{c}i\v{r}\'{\i}k, The first {R}obin eigenvalue with negative boundary parameter,
{\it Adv. Math.} \textbf{280} (2015), 322--339.

\bibitem[FL21]{FL18a}
P.~Freitas and R.~Laugesen,
From Neumann to Steklov and beyond, via Robin: the Weinberger way, \textit{ Amer. J. Math.} {\bf 143} (2021), 969--994.

\bibitem[FL20]{FL18b}
P. Freitas and R.\,S.~Laugesen,
From Steklov to Neumann and beyond, via Robin: the Szeg\H{o} way, 
\textit{Can. J. Math.}
{\bf 72} (2020), 1024--1043.



\bibitem[HKK01]{HKK01}
E.~Harrell, P.~Kr\"oger, and K.~Kurata, On the placement of an obstacle or a well so as to optimize the fundamental eigenvalue,
{\it SIAM J. Math. Anal.} \textbf{33} (2001), 240--259.

\bibitem[Hart64]{hart}
P.~Hartman, Geodesic parallel coordinates in the large, {\it Amer. J. Math.} \textbf{86} (1964), 705--727.



\bibitem[GL21]{GL19}
A.~Girouard and R.\,S.~Laugesen, Robin spectrum: two disks maximize the third eigenvalue,  {\it Indiana Univ. Math. J.} {\bf 70} (2021), 2711--2742.


\bibitem[K24]{K24}
E.~Krahn,
\"Uber eine von Rayleigh formulierte Minimaleigenschaft des Kreises, {\it Math. Ann.} \textbf{94} (1924), 97--100.

\bibitem[KL18]{KL17a}
D.~Krej\v{c}i\v{r}\'{\i}k and V.~Lotoreichik, Optimisation of the lowest {R}obin eigenvalue in the exterior of a compact set,
{\it J. Convex Anal.} {\bf 25} (2018), 319--337.

\bibitem[KL20]{KL17b}
D.~Krej\v{c}i\v{r}\'{i}k and V.~Lotoreichik, Optimisation of the lowest Robin eigenvalue in the exterior of a compact set, II: non-convex domains and higher dimensions, {\it Potential Anal.} {\bf 52} (2020), 601--614.
%


\bibitem[L19]{L19}
R.S.~Laugesen, The Robin Laplacian-spectral conjectures, rectangular theorems,
{\it J. Math. Phys.} {\bf 60} (2019), 121507, 32~pages.





\bibitem[R1877]{R77}
J. W. W. Rayleigh, {\it The  theory  of  sound},  Macmillan, London 1877.

\bibitem[V19]{V19}
A.~Vikulova, Parallel coordinates in three dimensions and sharp spectral isoperimetric inequalities, {\it to appear in Ric. Mat.}, \texttt{arXiv:1906.11141}.


\end{thebibliography}
